\documentclass{amsart}
\usepackage{amsfonts}
\usepackage{amssymb}
\usepackage{graphicx}
\usepackage{amsmath}

\newtheorem{theorem}{Theorem}

\newtheorem{corollary}[theorem]{Corollary}

\newtheorem{lemma}[theorem]{Lemma}

\input{tcilatex}

\begin{document}
\title{Leavitt Path Algebras of Finite Gelfand-Kirillov Dimension }
\author{Adel Alahmadi$^{(1)}$, Hamed Alsulami$^{(1)}$, S. K. Jain$^{(1,2)}$,
and Efim Zelmanov$^{(1,3)}$}
\address{1. Department of Mathematics\\
King Abdulaziz University\\
Jeddah, SA\\
2.Department of Mathematics\\
Ohio University, USA\\
3 Department of Mathematicsw\\
\ University of California, San Diego\\
Partially supported by the NSF}
\maketitle

\begin{abstract}
Groebner-Shirshov basis and Gelfand-Kirillov dimension of the Leavitt path
algebra are derived.
\end{abstract}

\section{Introduction.}

Leavitt path algebras were introduced in [AA] as algebraic analogs of graph
Cuntz-Kreiger C*-algebras. Since then they have received significant
attention from algebraists. In this paper we (i) find a Groebner-Shirshov
basis of a Leavitt path algebra, (ii) determine necessary and sufficient
conditions for polynomially bounded growth, and (iii) find Gelfand-Kirillov
dimension.

\section{\protect\bigskip Definitions and Terminologies}

A (directed) graph $\Gamma =(V,E,s,r)$ consists of two sets $V$ and $E,$
called vertices and edges respectively, and two maps $s,$ $r:E\rightarrow V$%
. The vertices $s(e)$ and $r(e)$ are referred to as the source and the range
of the edge $e$, respectively. The graph is called row-finite if for all
vertices $v\in V$, $|(s^{-1}(v))|<\infty .$ A\ vertex $v$ for which $%
(s^{-1}(v))$ is empty is called a sink. A path $p=e_{1}...e_{n}$ in a graph $%
\Gamma $ is a sequence of edges $e_{1}...e_{n}$ such that $%
r(e_{i})=s(e_{i+1})$ for $i=1,...,(n-1).$ In this case we say that the path $%
p$ starts at the vertex $s(e_{1})$ and ends at the vertex $r(e_{n}).$ If $%
s(e_{1})=$ $r(e_{n}),$ then the path is closed. If $p=e_{1}...e_{n}$ is a
closed path and the vertices $s(e_{1}),..,s(e_{n})$ are distinct, then the
subgraph $(s(e_{1}),...,$ $s(e_{n});e_{1},...,e_{n})$ of the graph $\Gamma $
is called a cycle.

Let $\Gamma $ be a row-finite graph and let $F$ be a field. The Leavitt path 
$F$-algebra $L(\Gamma )$ is the $F$-algebra presented by the set of
generators $\{v,$ $v\in V\},$ $\{e,$ $e^{\ast }|$ $e\in E\}$ and the set of
relators (1) $v_{i}v_{j}=\delta _{v_{i},v_{j}}v_{i}$ for all $v_{i},$ $%
v_{j}\in V;$ (2) $s(e)e=er(e)=e,$ $r(e)e^{\ast }=e^{\ast }s(e)=e^{\ast },$
for all $e\in E;$ (3) $e^{\ast }f=\delta _{e,f}r(e),$ for all $e,$ $f\in E;$
(4) $v=\sum_{s(e)=v}ee^{\ast },$ for an arbitrary vertex $v$ $\in
V\backslash \{sinks\}$.

The mapping which sends $v$ to $v,$\ for $v\in V,$ $e$ to $e^{\ast }$ and $%
e^{\ast }$ to $e,$ for $e\in E,$ extends to an involution of the algebra $%
L(\Gamma ).$ If $p=e_{1}...e_{n}$ is a path, then $p^{\ast }=e_{n}^{\ast
}...e_{1}^{\ast }.$

\section{ A Basis of $L(\Gamma )$\protect\bigskip}

For an arbitrary vertex $v$ which is not a sink, choose an edge $\gamma (v)$
such that $s(\gamma (v))=v.$ We will refer to this edge as special. In other
words, we fix a function $\gamma :V\backslash \{$sinks$\}\rightarrow E$ such
that $s(\gamma (v))=v$ for an arbitrary $v\in V\backslash \{$sinks$\}.$

\begin{theorem}
The following elements form a basis of the Leavitt path algebra $L(\Gamma ):$
\end{theorem}

(i) $v,$ $v\in V$, (ii) $p,$ $p^{\ast }$,where $p$ is a path in $\Gamma ,$
(iii) $pq^{\ast },$ where $p=e_{1}...e_{n}$ $,$ $q=f_{1}...f_{m}$ $,$ $%
e_{i}, $ $f_{j}\in E,$ are paths that end at the same vertex $%
r(e_{n})=r(f_{m}),$ with the condition that the last edges $e_{n}$ and $%
f_{m} $ are either distinct or equal, but not special.$\bigskip $

\begin{proof}
Recall that a well-ordering on a set is a total order (that is, any two
elements can be ordered) such that every non-empty subset of elements has a
least element.

As a first step, we will introduce a certain well-ordering on the set of
generators $X=V\cup E\cup E^{\ast }.$ Choose an arbitrary well-ordering on
the set of vertices $V.$ If $\ e,$ $f$ are edges and $s(e)<s(f)$ then $e<f.$
It remains to order edges that have the same source. Let $v$ be a vertex
which is not a sink. Let $e_{1},...,e_{k}$ be all the edges that originate
from $v.$ Suppose $e_{k}=\gamma (v).$ We order the edges as follows: $%
e_{1}<e_{2}<....<e_{k}=\gamma (v).$ Choose an arbitrary well-ordering on the
set $E^{\ast }.$ For arbitrary elements $v\in V,$ $e\in E,$ $f^{\ast }\in
E^{\ast },$ we let $v<e<f^{\ast }.$ Thus the set $X=V\cup E\cup E^{\ast }$
is well-ordered. Let $X^{\ast }$ be the set of all words in the alphabet $X.$
The length-lex order (see [B, Be]) makes $X^{\ast }$ a well-ordered set. For
all $v\in V$ and $e\in E$, we extend the set of relators (1) - (4) by (5): $%
ve=0,$ for $v\neq s(e);$ $ev=0,$ for $v\neq r(e);$ $ve^{\ast }=0,$ for $%
v\neq r(e);$ $e^{\ast }v=0,$ for $v\neq s(e).$ The straightforward
computations show that the set of relators (1) - (5) is closed with respect
to compositions (see [B, BE]). By the Composition-Diamond Lemma ([B, BE])
the set of irreducible words (not containing the leading monomials of
relators (1) - (5) as subwords) is a basis of $L(\Gamma ).$ This completes
the proof.\bigskip
\end{proof}

\section{Leavitt path algebras of polynomial growth}

Recall some general facts on the growth of algebras. Let $A$ be an algebra
(not necessarily unital), which is generated by a finite dimensional
subspace $V$. Let $V^{k}$ denote the span of all products $v_{1}\cdot \cdot
\cdot $ $v_{k},$ $v_{i}\in V,$ $k\leq n.$ Then $V=V^{1}\subset V^{2}\subset
\cdot \cdot \cdot ,$ $A=\cup _{n\geq 1}V^{n}$ and $g_{V(n)}=\dim
V^{n}<\infty .$ Given the functions $f,$ $g$ from $N$ $=\{1,2,...\}$ to the
positive real numbers $R_{+},$ we say that $f\preccurlyeq g$ if there exists 
$c\in N$ such that $f(n)$ $\leq cg(cn)$ for all $n$. If $f\preccurlyeq g$
and $g\preccurlyeq f$ then the functions $f,$ $g$ are said to be
asymptotically equivalent, and we write $f\thicksim g.$ If $W$ is another
finite dimensional subspace that generates $A$, then $g_{V(n)}$ $\thicksim $ 
$g_{W(n)}.$ If $g_{V(n)}$ is polynomially bounded, then we define the
Gelfand-Kirillov dimension of $A$ as $GKdim$ $A$ = $\lim \sup_{n\rightarrow
\infty }\frac{\ln g_{V(n)}}{\ln n}.$ The definition of $GK$-dimension does
not depend on a choice of the generating space $V$ as long as $dim$ $V$ 
\TEXTsymbol{<} $\infty .$ If the growth of $A$ is not polynomially bounded,
then $GKdim$ $A=$ $\infty .$

We now focus on finitely generated algebras and we will assume that the
graph $\Gamma $ is finite. Let $C_{1},$ $C_{2}$ be distinct cycles such that 
$V(C_{1})\cap V(C_{2})$ $\neq $ $\phi .$ Then we can renumber the vertices
so that $C_{1}=(v_{1},\cdot \cdot \cdot ,v_{m};e_{1},\cdot \cdot \cdot
,e_{m}),$ $C_{2}=(w_{1},\cdot \cdot \cdot ,w_{n};f_{1},\cdot \cdot \cdot
,f_{n}),$ $v_{1}=w_{1}.$ Let $p=e_{1}\cdot \cdot \cdot e_{m},$ and $%
q=f_{1}\cdot \cdot \cdot f_{n}.$

\begin{lemma}
The elements $p,$ $q$ generate a free subalgebra in $L(\Gamma )$.

\begin{proof}
By Theorem 1, different paths viewed as elements of $L(\Gamma )$ are
linearly independent. If $u_{1},$ $u_{2}$ are different words in two
variables, then $u_{1}(p,q)$ and $u_{2}(p,q)$ are different paths. Indeed,
cutting out a possible common beginning, we can assume that $u_{1},$ $u_{2}$
start with different letters, say, $u_{1}(p,q)=p\cdot \cdot \cdot ,$ $%
u_{2}(p,q)=q\cdot \cdot \cdot .$ If $m>n$ then the path $u_{2}(p,q)$ returns
to the vertex $v_{1}$ at the $n$-th step, whereas $u_{1}(p,q)$ does not. If $%
m=n,$ then the left segments of length $m$ of $u_{1}(p,q),$ $u_{2}(p,q)$ are
different. This proves the lemma.
\end{proof}
\end{lemma}

\begin{corollary}
If two distinct cycles have a common vertex, then $L(\Gamma )$ has
exponential growth.
\end{corollary}

From now on we will assume that any two distinct cycles of the graph $\Gamma 
$ do not have a common vertex.

For two cycles $C^{\prime },$ $C^{^{\prime \prime }},$ we write $C^{\prime
}\Longrightarrow C^{^{\prime \prime }},$ if there exists a path that starts
in $C^{\prime }$and ends in $C^{^{\prime \prime }}$.

\begin{lemma}
If $C^{\prime },$ $C^{^{\prime \prime }}$are two cycles such that $C^{\prime
}\Longrightarrow $ $C^{^{\prime \prime }},$ and $C^{^{\prime \prime
}}\Longrightarrow $ $C^{\prime },$ then $C^{\prime }=C^{^{\prime \prime }}$.

\begin{proof}
Choose a path $p$ that starts in $C^{\prime }$ and ends in $C^{^{\prime
\prime }}$. Similarly, choose a path $q$ that starts in $C^{^{\prime \prime
}}$ and finishes in $C^{\prime }$. There exists also a path $p^{\prime }$ on 
$C^{^{\prime \prime }}$, which connects $r(p)$ with $s(q)$ and a path $%
q^{\prime }$ on $C^{\prime }$, which connects $r(q)$ with $s(p)$. Now, $%
pp^{\prime }qq^{\prime }$ is a closed path, which visits both $C^{\prime }$
and $C^{^{\prime \prime }}$. Let $t$ be a closed path with this property
(visiting both $C^{\prime }$ and $C^{^{\prime \prime }}$) having a minimal
length. Write $t=$ $e_{1}\cdot \cdot \cdot e_{n},$ $e_{i}\in E.$ We claim
that the vertices $s(e_{1}),\cdot \cdot \cdot ,s(e_{n})$ are all distinct,
thus $t=(s(e_{1}),\cdot \cdot \cdot ,s(e_{n});e_{1},\cdot \cdot \cdot
,e_{n}) $ is a cycle. Assuming the contrary, let $s(e_{i})=s(e_{j}),$ $1\leq
i<j\leq n,$ and $j-i$ is minimal with this property. Then $t^{\prime }=$ $%
(s(e_{i}),s(e_{i+1}),\cdot \cdot \cdot ,s(e_{j});e_{i},e_{i+1},\cdot \cdot
\cdot ,e_{j-1})$ is a cycle. Let us "cut it out", that is, consider the path 
$t^{^{\prime \prime }}=e_{1}\cdot \cdot \cdot $ $e_{i-1}e_{j}\cdot \cdot
\cdot $ $e_{n}.$ This path is shorter than $t$. Hence $t^{^{\prime \prime }}$%
can not visit both $C^{\prime }$ and $C^{^{\prime \prime }}$. Suppose that $%
t^{^{\prime \prime }}$ does not visit $C^{\prime }$. Then at least one of
the vertices $s(e_{i}),\cdot \cdot \cdot ,s(e_{j-1})$ lies in $C^{\prime }$.
Since two intersecting cycles coincide, it implies that $t^{\prime
}=C^{\prime }$, hence $s(e_{j})$ lies in $C^{\prime }$. This contradicts our
assumption that $t^{^{\prime \prime }}$ does not visit $C^{\prime }$. Hence $%
t=C^{\prime }=C^{^{\prime \prime }}$. This proves the lemma.
\end{proof}
\end{lemma}

A sequence of distinct cycles $C_{1},\cdot \cdot \cdot ,$ $C_{k}$ is a chain
of length $k$ if $C_{1}\Longrightarrow \cdot \cdot \cdot \Longrightarrow $ $%
C_{k}.$ The chain is said to have an exit if the cycle $C_{k}$ has an exit
(see [AA]), that is, if there exists an edge $e$ such that $s(e)\in
V(C_{k}), $ but $e$ does not belong to $C_{k}.$Let $d_{1}$ be the maximal
length of a chain of cycles in $\Gamma $, and let $d_{2}$ be the maximal
length of chain of cycles with an exit. Clearly, $d_{2}\leq d_{1}.$

\begin{theorem}
Let $\Gamma $ be a finite graph.
\end{theorem}

(1) The Leavitt path algebra $L(\Gamma )$ has polynomially bounded growth if
and only if any two distinct cycles of $\ \Gamma $ do not have a common
vertex;

(2) If $d_{1}$ is the maximal length of a chain of cycles in $\Gamma $, and $%
d_{2}$ is the maximal length of chain of cycles with an exit, then $GK$ $dim$
$L(\Gamma )$ = $max(2d_{1}-1,2d_{2})$.

\begin{proof}
As in the proof of Theorem 1 we consider the generating set $X=V\cup E\cup
E^{\ast }$ of $L(\Gamma )$. Let $E^{\prime }$ be the set of edges that do
not belong to any cycle. Let $P^{\prime }$ be the set of all paths that are
composed from edges from $E^{^{\prime }}$. Then an arbitrary path from $%
P^{\prime }$ never arrives to the same vertex twice. Hence, $|P^{\prime }|$ $%
<\infty .$

By Theorem 1 the space $Span(X^{n})$ is spanned by elements of the following
types:

(1) a vertex,

(2) a path $p=p_{1}^{\prime }p_{1}p_{2}^{\prime }p_{2}\cdot \cdot \cdot
p_{k}p_{k+1}^{\prime },$ where $p_{i}$ is a path on a cycle $C_{i},$ $1\leq
i\leq k,$ $C_{1}\Longrightarrow \cdot \cdot \cdot \Longrightarrow $ $C_{k}$
is a chain, $p_{i}^{\prime }\in P^{\prime },$ $length(p)$ $\leq n,$

(3) $p^{\ast },$ where $p$ is a path of the type (2),

(4) $pq^{\ast },$ where $p=p_{1}^{\prime }p_{1}p_{2}^{\prime }p_{2}\cdot
\cdot \cdot p_{k}p_{k+1}^{\prime },$ $q=q_{1}^{\prime }q_{1}q_{2}^{\prime
}\cdot \cdot \cdot q_{s}q_{s+1}^{\prime };$ $p_{i},$ $q_{j}$ are paths on
cycles $C_{i},$ $D_{j\text{ }}$ respectively and $C_{1}\Longrightarrow \cdot
\cdot \cdot \Longrightarrow $ $C_{k},$ $D_{1}\Longrightarrow \cdot \cdot
\cdot \Longrightarrow $ $Ds$ are chains ; $p_{i}^{\prime },q_{j}^{\prime
}\in P^{\prime },$ $length(p)+length(q)\leq n$ with $r(p)=r(q).$ We will
further subdivide this case into two subcases:

(4.1)\ $r(p)\notin V(C_{k})\cup V(D_{s});$\ 

(4.2) $r(p)\in V(C_{k})\cup V(D_{s}).$

We will estimate the number of products of $length$ $\leq n$ in each of the
above cases and then use the following elementary fact:

Let $(a_{n})_{n\in N}$ \ be the sum of $s$ sequences $(a_{in})_{n\in N},$ $%
1\leq i\leq s,$ $a_{in}>0.$ Then

lim sup$_{n\rightarrow \infty }$ $\frac{\ln \text{ }a_{n}}{\ln \text{ }n}%
=\max (\lim $ $\sup_{n\rightarrow \infty }\frac{\ln \text{ }a_{in}}{\ln 
\text{ }n},$ $1\leq i\leq s)$

Let us estimate the number of paths of the type (2). Fix a chain $C_{1}$ $%
\Longrightarrow $ $C_{2}$ $\Longrightarrow ...\Longrightarrow C_{k}.$ If $%
C=(v_{1},....,v_{m};e_{1},...,e_{m})$ is a cycle, let $P_{C}=e_{1}...e_{m}.$
For a given cycle there are $m$ such paths depending upon the choice of the
starting point $v_{1}.$

Let $|(V(C_{i})|=m_{i}$ and let $P_{C_{i}}$ be any one of the $m_{i}$ paths
described above. Then an arbitrary path on $C_{i}$ can be represented as $%
u^{\prime }P_{i}^{l}u^{\prime \prime },$ where $length$ $(u^{\prime })$, $%
length$ $(u^{^{\prime \prime }})$ $\leq $ $m_{i}-1.$ Hence every path of the
type (2) which corresponds to the chain $C_{1}$ $\Longrightarrow $ $C_{2}$ $%
\Longrightarrow ...\Longrightarrow $ $C_{k}$ can be represented as $%
p_{1}^{\prime }u_{1}^{\prime }P_{C_{1}}^{l_{1}}u_{1}^{\prime \prime
}...p_{k}^{\prime }u_{k}^{\prime }P_{Ck}^{l_{k}}u_{k}^{\prime \prime
}p_{k+1}^{\prime },$ where $p_{i}^{\prime }\in P_{i}^{\prime }$ and $length$ 
$(u_{i}^{\prime }),$ $length$ $(u_{i}^{\prime \prime })\leq m_{i}-1.$
Clearly, $m_{1}l_{1}+...+m_{k}l_{k}\leq n.$ This implies that the number of
such paths $\preccurlyeq $ $n^{k}\leq n^{d_{1}}.$ On the other hand,
choosing a chain $C_{1}$ $\Longrightarrow $ $C_{2}$ $\Longrightarrow
...\Longrightarrow $ $C_{d_{1}}$ of length $d_{1}$, we can construct $\sim $ 
$n^{d_{1}}$ paths\ of length $\leq n.$ The case (3) is similar to the case
(2).

Consider now the elements of length $\leq n$ of the type $pq^{\ast },$ $%
r(p)=r(q);$ the path $p$ passes through the cycles of the chain $C_{1}$ $%
\Longrightarrow $ $C_{2}$ $\Longrightarrow ...\Longrightarrow $ $C_{k}$ on
the way, the path $q$ passes through the cycles of the chain $D_{1}$ $%
\Longrightarrow $ $D_{2}$ $\Longrightarrow ...\Longrightarrow $ $D_{s}$ on
the way and so $p=p_{1}^{\prime }p_{1}p_{2}^{\prime }...p_{k}p_{k+1}^{\prime
},$ where \ $p_{i}^{\prime }\in P_{i}^{\prime },$ each $p_{i}$ is a path on
the cycle $C_{i\text{ }}.$ Similarly, $q=q_{1}^{\prime }q_{1}q_{2}^{\prime
}...q_{s}q_{s+1}^{\prime }.$ Arguing as above, we see that for fixed chains $%
C_{1}$ $\Longrightarrow $ $C_{2}$ $\Longrightarrow ...\Longrightarrow $ $%
C_{k}$ and $D_{1}$ $\Longrightarrow $ $D_{2}$ $\Longrightarrow
...\Longrightarrow $ $D_{s}$ , the number of such paths $\preccurlyeq
n^{k+s}.$

Suppose that the vertex $v=r(p)=r(q)$ does not lie in $V(C_{k})\cup V(D_{s})$
. Then both cycles $C_{k}$ and $D_{s}$ have exits. Hence the number of paths
of type (4.1) is $\leq n^{2d_{2}\text{ }}.$ On the other hand, let $C_{1}$ $%
\Longrightarrow $ $C_{2}$ $\Longrightarrow ...\Longrightarrow $ $C_{d_{2}}$
be a chain and let $e$ be an exit of the cycle $C_{d_{2}}.$ Select\ paths $%
p_{2}^{\prime }$, ...,$p_{d_{2}}^{\prime }$, where $p_{i}^{\prime }$
connects $C_{i-1}$ to $C_{i}$ , $p_{i}^{\prime }\in P^{\prime }.$

Select a path $u_{1}^{\prime \prime }$ on the cycle $C_{1\text{ }}$ which
connects $r(P_{C_{1}}$ $)$ to $s(p_{2}^{\prime })$ , a path $u_{2}^{\prime }$
in $C_{2}$which connects $r(p_{2}^{\prime })$ to $s(P_{C_{2}}),$ a path $%
u_{2}^{\prime \prime }$ on $C_{2}$ which connects $r(P_{C_{2}})$ to $%
s(P_{3}^{\prime })$, and so on. The path $u_{d_{2}}^{\prime \prime }$
connects $r(P_{C_{d_{2}}})$ \ to $s(e).$

Among the edges from $s^{-1}(s(e))$ choose a special one $\gamma (s(e))$
different from $e.$ Then by Theorem 1, the elements

$P_{C_{1}}^{l_{1}}u_{1}^{\prime \prime }p_{2}^{\prime }u_{2}^{\prime
}P_{C2}^{l_{2}}u_{2}^{\prime \prime }p_{3}^{\prime
}...P_{C_{d_{2}}}^{l_{d_{2}}}u_{d_{3}}^{\prime \prime }ee^{\ast
}(u_{2}^{\prime \prime })^{\ast }(P_{C_{d_{2}}^{\ast
}})^{l_{d_{2}+1}}...(u_{1}^{\prime \prime })^{\ast }(P_{C_{1}})^{l_{2d_{2}}}$
$,$ $l_{i}\geq 1,$ \ \ (A),

are linearly independent. Let $m$ be the total length of all elements other
than $P_{C\ _{i}}^{l_{i}},$ ($P_{C_{i}^{\ast }})^{l_{2d_{2}-i+1}}$. The
number of elements in the inequality (A) above is the number of nonnegative
integral solutions of the inequality

$\sum_{i=1}^{d_{2}}m_{i}(l_{i}+l_{2d_{2}-i+1})\leq n-m,$ which is $\thicksim
n^{2d_{2}}.$\bigskip

Now suppose that the vertex $v=r(p)=r(q)$ lies in $C_{k}.$ Assume \ at first
that $C_{k}\neq D_{s}.$ Then the chain $D_{1}$ $\Longrightarrow $ $D_{2}$ $%
\Longrightarrow ...\Longrightarrow $ $D_{s}$ has an exit. If $k\leq s,$ then
the number of the paths of this type is $\leq n^{k+s}\leq n^{2d_{2}}.$

If $s<k,$ then $n^{k+s}\leq n^{2k-1}\leq n^{2d_{1}}-1.$

Next, let $C_{k}=D_{s}.$ It means that the paths $p_{k+1}^{\prime
},q_{s+1}^{\prime }$ are empty; $p_{k}$ and $q_{s}$ are both paths on the
cycle $C_{k}$ and in this case we have,

$(i)$ $p_{k}q_{s}^{\ast }$ = $u,$ \ if $p_{k}=uq_{s}$, is a path on $C_{k},$

$(ii)p_{k}q_{s}^{\ast }=$ $u^{\ast },$ \ if $q_{s}=up_{k}$, is a path on $%
C_{k}$\ ,and

$(iii)$\ $p_{k}q_{s}^{\ast }=$\ $0,$\ \ otherwise.

The number of such elements $pq^{\ast }$ is $\preccurlyeq n^{k+s-1}\leq
n^{2d_{1}}-1.$

On the other hand, let $C_{1}$ $\Longrightarrow $ $C_{2}$ $\Longrightarrow
...\Longrightarrow $ $C_{d_{2}}$ be a chain of cycles. Select paths $%
p_{2}^{\prime },...,p_{d_{1}}^{\prime }\in P^{\prime },$ $p_{i}^{\prime }$
connects $C_{i-1}$ to $C_{i};$ $u_{i}^{\prime }$, $u_{i}^{\prime \prime }$
are paths on the cycle $C_{i}$ such that $P_{C_{1}}u_{1}^{\prime \prime }$ $%
p_{2}^{\prime }u_{2}^{\prime }P_{C_{2}}u_{2}^{\prime \prime }p_{3}^{\prime
}...P_{C_{d_{1}}}\neq $ $0.$ By Theorem 1, the elements $%
P_{C_{1}}^{l_{1}}u_{1}^{\prime \prime }p_{2}^{\prime }u_{2}^{\prime
}P_{C_{2}}^{l_{2}}u_{2}^{\prime \prime }p_{3}^{\prime }...u_{d_{1}}^{\prime
}P_{C_{d_{1}}}^{l_{d_{1}}}(u_{d_{1}}^{\prime })^{\ast
}(P_{C_{d_{1}-1}}^{\ast })^{l_{d_{1}+1}}...(P_{C_{1}}^{\ast
})^{l_{2d_{1}-1}} $ are linearly independent provided that $l_{i}$ $\geq $ $%
1 $, $1\leq i$ $\leq 2d_{1}-1.$ The number of these elements is $\sim $\ $%
n^{2d_{1}-1.}$ This proves Theorem 2.
\end{proof}

References.

[AA] G. Abrams, G. Aranda Pino, The Leavitt path algebra of a graph, J.
Algebra 293 (2005), 319-334

[B] L.A. Bokut, Imbeddings in simple associative algebras, Algebra i Logica,
15 (1976), 117-142.

[Be] G.M. Bergman, The Diamond Lemma for ring theory, Adv. in Math. 29
(1978), 178-218.

[KL] G.R.Krause, T.H. Lenagan, Growth of Algebras and Gelfand-Kinillow
Dimension, revised edition. Graduate studies in Mathematics, 22, AMS,
Providence, RI, 2000

\end{document}